\documentclass[a4paper,11pt,reqno]{amsart}
\usepackage[utf8]{inputenc}
\usepackage[T1]{fontenc}
\usepackage{latexsym,amsmath,amsfonts,amscd,amssymb,amsthm,mathrsfs}
\usepackage{pifont}
\usepackage{eucal}
\usepackage{enumerate}
\usepackage[all]{xy}
\usepackage{color}
\usepackage{tikz}
\usepackage{listings}
\usepackage{fullpage}

\definecolor{cobalt}{RGB}{61,89,171}
\usepackage[colorlinks,citecolor=cobalt,linkcolor=cobalt,urlcolor=cobalt,pdfpagemode=UseNone,backref = page]{hyperref}

\newcommand{\tens}[1]{
  \mathbin{\mathop{\otimes}\displaylimits_{#1}}
}

\newcommand{\etalchar}[1]{$^{#1}$}

\theoremstyle{plain}
\newtheorem{theorem}{Theorem}[section]
\newtheorem{proposition}[theorem]{Proposition}

\theoremstyle{definition}

\theoremstyle{remark}

\begin{document}

\title{SageMath experiments in Differential and Complex Geometry}

\author[D. Angella]{Daniele Angella}

\address[D. Angella]{
Dipartimento di Matematica e Informatica "Ulisse Dini"\\
Università di Firenze\\
viale Morgagni 67/a\\
50134 Firenze\\
Italy\\
}

\email{daniele.angella@gmail.com}
\email{daniele.angella@unifi.it}

\urladdr{http://sites.google.com/site/danieleangella/}

\thanks{The author is supported by the Project FIRB ``Geometria Differenziale e Teoria Geometrica delle Funzioni'', by the Project SIR 2014 AnHyC ``Analytic aspects in complex and hypercomplex geometry'' (code RBSI14DYEB), and by GNSAGA of INdAM}

\date{\today}

\begin{abstract}
This note summarizes the talk by the author at the workshop ``Geometry and Computer Science'' held in Pescara in February 2017.
We present how SageMath can help in research in Complex and Differential Geometry, with two simple applications, which are not intended to be original.
We consider two "classification problems" on quotients of Lie groups, namely, "computing cohomological invariants" \cite{angella-franzini-rossi, latorre-ugarte-villacampa}, and "classifying special geometric structures" \cite{angella-bazzoni-parton}, and we set the problems to be solved with SageMath \cite{sage}.
\end{abstract}

\maketitle

\section*{Introduction}

Complex Geometry is the study of manifolds locally modelled on the linear complex space $\mathbb{C}^n$. A natural way to construct compact complex manifolds is to study the projective geometry of $\mathbb{C}\mathbb{P}^n=\mathbb{C}^{n+1}\setminus\{0\}\slash \mathbb{C}^\times$. In fact, analytic submanifolds of $\mathbb{C}\mathbb{P}^n$ are equivalent to algebraic submanifolds \cite{serre-gaga}. On the one side, this means that both algebraic and analytic techniques are available for their study. On the other side, this means also that this class of manifolds is quite restrictive. In particular, they do not suffice to describe some Theoretical Physics models {\itshape e.g.} \cite{strominger}. Since \cite{thurston}, new examples of complex non-projective, even non-K\"ahler manifolds have been investigated, and many different constructions have been proposed. One would study classes of manifolds whose geometry is, in some sense, combinatorically or algebraically described, in order to perform explicit computations. In this sense, great interest has been deserved to homogeneous spaces of nilpotent Lie groups, say {\em nilmanifolds}, whose geometry \cite{belgun, fino-grantcharov} and cohomology \cite{nomizu} is encoded in the Lie algebra. In non-K\"ahler geometry, they often have a role of toy-models: to prove or disprove conjecture {\itshape e.g.} \cite{fernandez-munoz}, to explicitly solve specific equations {\itshape e.g.} \cite{buzano-fino-vezzoni, tosatti-weinkove-jussieu}, to provide examples of models {\itshape e.g.} \cite{fernandez-ivanov-ugarte-villacampa, ugarte-villacampa-asian, ugarte-villacampa-forum, otal-ugarte-villacampa}, to get an idea of possible general results. Algebraic techniques have allowed to classify the linear models for $6$-dimensional nilmanifolds with invariant complex structures \cite{salamon}, and the invariant complex structures on them \cite{ceballos-otal-ugarte-villacampa}, so paving the way for a complete description of cohomological \cite{latorre-ugarte-villacampa, angella-franzini-rossi} and metric properties.

In this note, we propose two problems at a research level. They appeared in \cite{angella-franzini-rossi}, respectively \cite{angella-bazzoni-parton}, and they are concerned with the description of cohomological and metric properties of nilmanifolds and Lie algebras. In solving them, we have made use of symbolic computations, especially using SageMath \cite{sage} at SageMathCloud. SageMath is a "free, open-source math software that supports research and teaching in algebra, geometry, number theory, cryptography, numerical computation, and related areas" \cite{sage-tutorial}; see also \cite{bard}. In this spirit, we propose here our simple solutions, with the aim to serve as a basic introduction for researchers in Complex Geometry. In fact, the algorithms here proposed are {\em not}, in any way, complete or efficient.

\bigskip

\noindent{\sl Acknowledgments.}
This note has been written for the Workshop ``Geometry and Computer Science'' held in Pescara on February 8th--10th, 2017, \url{http://www.sci.unich.it/gncs2017/}.
The author warmly thanks the organizers for the kind invitation and hospitality, and also all the participants for the environment they contributed to.
During the years, the author had the opportunity to learn how to use symbolic mathematics software systems thanks to Nicola Enrietti, Antonio Otal, Federico A. Rossi, Miguel Angel Marco Buzunariz, among others. Thanks also to Giovanni Bazzoni, Maria Giovanna Franzini, Maurizio Parton, Federico A. Rossi, for the opportunity to collaborate with them on the projects \cite{angella-franzini-rossi, angella-bazzoni-parton}.

\section{Preliminaries on homogeneous spaces of Lie groups and Lie algebras}
Our primary objects of study are {\em nilmanifolds} $\left.\Gamma\backslash G \right.$, namely, compact quotients of connected simply-connected nilpotent Lie groups $G$ by co-compact discrete subgroups $\Gamma$. We recall that a {\em Lie group} is a group endowed with a structure of differentiable manifold such that the group operations are smooth. In fact, we are interested in nilmanifolds also because many of their properties can be reduced at the linear level of the corresponding Lie algebras; see {\itshape e.g.} the averaging trick in \cite[Theorem 7]{belgun}. We recall that the {\em Lie algebra} $\mathfrak{g}$ of a Lie group $G$ is the tangent space at the identity $e \in G$. Therefore, it is a vector space endowed with a Lie bracket $[\_,\_]\colon \mathfrak{g}\times\mathfrak{g}\to\mathfrak{g}$, that is, a bilinear skew-symmetric map satisfying the Jacobi identity: $[x,[y,z]]+[y,[z,x]]+[z,[x,y]]=0$ for any $x,y,z\in\mathfrak{g}$. Focusing on the cotangent space instead of the tangent space, namely, on the dual $\mathfrak{g}^\vee$ instead of $\mathfrak{g}$, we can construct a linear operator $d\colon \mathfrak{g}^\vee \to \wedge^2\mathfrak{g}^\vee$ by
$$ (d \alpha)(x,y) := -\alpha([x,y]) .$$
We can extend it to a differential operator on the complex $\wedge^\bullet \mathfrak{g}^\vee$. We do this by requiring the Leibniz rule: $d(\alpha\wedge\beta)=d(\alpha)\wedge\beta+(-1)^a\alpha\wedge d(\beta)$, for $\alpha\in\wedge^a\mathfrak{g}^\vee$ and $\beta\in\wedge^b\mathfrak{g}^\vee$. Namely, $d$ is a linear operator of degree $1$, that is, $d \colon \wedge^{\bullet}\mathfrak{g}^\vee \to \wedge^{\bullet+1}\mathfrak{g}^\vee$, satisfying the Leibniz rule, and with the property $d^2=0$. The last properties allows to consider the vector space
$$ H^\bullet(\mathfrak{g};\mathbb{R}) := \frac{\ker d}{\mathrm{im}\, d} .$$
We notice that any element in $\wedge^\bullet \mathfrak{g}^\vee$ can be thought of as a linear form on $\mathfrak{g}=T_eG$. By translating it with left-multiplication on $G$, we recover a differential form on $G$, which is invariant under the action of left-multiplication of $G$ on itself. So we get an inclusion $\wedge^\bullet \mathfrak{g}^\vee \to \wedge^\bullet T^\vee G$, whence also $\wedge^\bullet\mathfrak{g}^\vee \to \wedge^\bullet X$.
The following theorem by K. Nomizu allows us to compute some topological invariant of $X$ at the linear level of $\mathfrak{g}$. This makes the computation of the Betti numbers of $X$ just a matter of linear algebra. Recall that, in general, on compact differentiable manifolds, we can make use of Hodge theory to reduce the computation of the Betti numbers to a system of partial differential equations.

\begin{theorem}[{\cite[Theorem 1]{nomizu}}]
Let $X=\left.\Gamma\backslash G \right.$ be a nilmanifold, with associated Lie algebra $\mathfrak{g}$. Then the inclusion $\wedge^\bullet\mathfrak{g}^\vee\to \wedge^\bullet X$ induces the isomorphism $H^\bullet(\mathfrak{g};\mathbb{R}) \to H^\bullet(X;\mathbb{R})$, where $H^\bullet(X;\mathbb{R})$ denotes the de Rham cohomology of $X$. In particular, the Betti numbers of $X$ are (varying $j\in\mathbb{Z}$)
$$ b_j = \dim H^j(\mathfrak{g};\mathbb{R}) . $$
\end{theorem}

We consider further structures on nilmanifolds, respectively Lie algebras. As a general pattern: we consider an invariant structure on $X$; (if not invariant, we can possibly symmetrize it to an invariant one \cite[Theorem 7]{belgun};) here, invariant means that the structure is encoded in a linear structure on the Lie algebra; we associate some cohomological invariant to that structure; we extend the Nomizu theorem to such cohomologies.

For example, let us consider an {\em invariant complex structure} on $X$. This means a linear complex structure on $\mathfrak{g}$, that is, a structure of complex vector space. In other words, we have an endomorphism $J\in\mathrm{End}(\mathfrak{g})$ such that $J^2=-\mathrm{id}$. We can extend it to a linear complex structure $J$ on $\mathfrak{g}^\vee$ by $J\alpha:=\alpha(J^{-1}\_)$. This extends also to a morphism $J$ of the exterior algebra $\wedge^\bullet \mathfrak{g}^\vee$, by $J\phi:=\phi(J^{-1}\_, \dots, J^{-1}\_)$.
By complexifying $\mathfrak{g}$, we have a splitting into eigenspaces for $J$, namely, $\mathfrak{g}\tens{\mathbb{R}}\mathbb{C}=\mathfrak{g}^{1,0}\oplus\mathfrak{g}^{0,1}$ where $J\lfloor_{\mathfrak{g}^{1,0}}=\sqrt{-1}$ and $J\lfloor_{\mathfrak{g}^{0,1}}=-\sqrt{-1}$. Analogously, we can split the complexified exterior algebra $\wedge^\bullet \mathfrak{g}^\vee\tens{\mathbb{R}}\mathbb{C}=\bigoplus_{p+q=\bullet} \wedge^p(\mathfrak{g}^{1,0})^\vee \otimes \wedge^q(\mathfrak{g}^{0,1})^\vee=:\bigoplus_{p+q=\bullet} \wedge^{p,q}\mathfrak{g}^\vee$. Also the differential split into $d=\partial+\overline\partial$, where $\partial\colon \wedge^{\bullet,\bullet}\mathfrak{g}^\vee \to \wedge^{\bullet+1,\bullet}\mathfrak{g}^\vee$ and $\overline\partial \colon \wedge^{\bullet,\bullet}\mathfrak{g}^\vee \to \wedge^{\bullet,\bullet+1}\mathfrak{g}^\vee$. The condition $d^2=0$ yields $\partial^2=\overline\partial^2=\partial\overline\partial+\overline\partial\partial=0$. So we can consider the cohomology
$$ H^{\bullet,\bullet}_{\overline\partial}(\mathfrak{g}):=\frac{\ker\overline\partial}{\mathrm{im}\, \overline\partial}. $$
An analogue of the Nomizu theorem in the complex setting allows to compute some holomorphic invariants of $X$ thanks to $H^{\bullet,\bullet}(\mathfrak{g})$ \cite{console, rollenske-survey}.

\begin{theorem}[{\cite[Theorem 1]{sakane}, \cite[Main Theorem]{cordero-fernandez-gray-ugarte-2000}, \cite[Theorem 2, Remark 4]{console-fino}, \cite[Theorem 1.10]{rollenske}}]\label{thm:dolbeault-invariant}
Let $X=\left. \Gamma \middle\backslash G \right.$ be a nilmanifold, with associated Lie algebra $\mathfrak{g}$. Let $J$ be an invariant complex structure on $X$.
Assume some suitable conditions on the complex structure ({\itshape e.g.} holomorphically parallelizable; Abelian; nilpotent; rational).
Then the inclusion $\wedge^{\bullet,\bullet}\mathfrak{g}^\vee \to \wedge^{\bullet,\bullet}X$ induces the isomorphism $H^{\bullet,\bullet}_{\overline\partial}(\mathfrak{g})\to H^{\bullet,\bullet}_{\overline\partial}(X)$, where $H^{\bullet,\bullet}_{\overline\partial}(X)$ denotes the Dolbeault cohomology of $X$. In particular, the Hodge numbers of $X$ are (varying $p,q\in\mathbb{Z}$):
$$ h^{p,q}=\dim H^{p,q}_{\overline\partial}(\mathfrak{g}) .$$
\end{theorem}

Further cohomological invariants can be defined, and the same result as Nomizu's and Sakane's theorems applies \cite{angella-1}: the {\em Bott-Chern cohomology} and the {\em Aeppli cohomology} (either of a complex manifold, or of a Lie algebra with linear complex structure) are defined as
$$ H^{\bullet,\bullet}_{BC} := \frac{\ker\partial\cap\ker\overline\partial}{\mathrm{im}\,\partial\overline\partial} , \qquad H^{\bullet,\bullet}_{A} := \frac{\ker\partial\overline\partial}{\mathrm{im}\,\partial+\mathrm{im}\,\overline\partial} . $$

As another example, we consider {\em symplectic} and {\em locally conformally symplectic (lcs)} structures on nilmanifolds. Any such structure, thanks to Belgun symmetrization \cite[Theorem 7]{belgun}, (and thanks to \cite[Theorem 1]{nomizu},) gives an invariant structure, and so a linear structure on the Lie algebra. Recall that a symplectic structure is a non-degenerate $2$-form being $d$-closed. Linearly, it is given by $\omega\in\wedge^2\mathfrak{g}^\vee$ such that $\omega^{\frac{1}{2}\dim\mathfrak{g}}\neq0$ and $d\omega=0$. A locally conformally symplectic structure is a structure that is locally conformally to a symplectic structure, that is, a non-degenerate $2$-form satisfying $d\omega=\vartheta\wedge\omega$ for some $d$-closed $1$-form $\vartheta$, which is called the {\em Lee form}. In this case, we can define the differential $d_{\vartheta}:=d-\vartheta\wedge\_$, and its cohomology
$$ H^\bullet_{\vartheta}(\mathfrak{g}):=\frac{\ker d_\vartheta}{\mathrm{im}\,d_\vartheta} .$$
We have the Hattori theorem.

\begin{theorem}[{\cite[Corollary 4.2]{hattori}}]
Let $X=\left. \Gamma \middle\backslash G \right.$ be a nilmanifold, with associated Lie algebra $\mathfrak{g}$. Let $\omega$ be an invariant locally conformally structure on $X$, with Lee form $\vartheta$.
Then the inclusion $\wedge^{\bullet}\mathfrak{g}^\vee \to \wedge^{\bullet}X$ induces the isomorphism $H^{\bullet}_{d_{\vartheta}}(\mathfrak{g})\to H^{\bullet}_{d_{\vartheta}}(X)$, where $H^{\bullet}_{d_{\vartheta}}(X)$ denotes the Morse-Novikov cohomology of $X$. In particular, the Novikov Betti numbers of $X$ are (varying $k\in\mathbb{Z}$):
$$ b_{k}^{\text{Nov}} = \dim H^{k}_{d_\vartheta}(\mathfrak{g}) .$$
\end{theorem}

\section{Preliminaries on Clifford Algebras module in SageMath}\label{sec:preliminaries-sage}
We will use the class \lstinline{ExteriorAlgebra} and the class \lstinline{ExteriorAlgebraCoboundary} in the module \lstinline{Clifford Algebras}; see the corresponding chapter in the Sage Reference manual \cite{sage-reference} for more details.

\medskip

In particular, \lstinline{E = ExteriorAlgebra(R,'e',n)} creates an exterior algebra over the free module of rank \lstinline{n} over the ring \lstinline{R}, referred by \lstinline{E}. Possible rings are $\lstinline{RR}$ for $\mathbb{R}$; $\lstinline{QQ}$ for $\mathbb{Q}$; $\lstinline{SR}$ for the symbolic ring. The generators of the algebra,
\begin{lstlisting}
sage: E.gens()
\end{lstlisting}
are named $\lstinline{e}_0$, $\lstinline{e}_1$, \dots, $\lstinline{e}_{\lstinline{n}}$, and they are subject to the relations $\lstinline{e}_j\wedge\lstinline{e}_k=-\lstinline{e}_k\wedge\lstinline{e}_j$. Here, "$\wedge$" is the exterior product, which is denoted by "\lstinline{*}".
The elements of \lstinline{E} are of type $\lstinline{e}_{j_1} \wedge \cdots \wedge \lstinline{e}_{j_k}$; for example, \lstinline{E.basis(p)} returns the basis of $\wedge^{\lstinline{p}} \lstinline{R}\langle \lstinline{e}_0, \dots, \lstinline{e}_{\lstinline{n}}\rangle$ as \lstinline{R}-module. The elements of an exterior algebra are instances of the \lstinline{Element} class. The method \lstinline{interior_product(x)} returns the interior product of \lstinline{self} with \lstinline{x}. The method \lstinline{hodge_dual()} returns the Hodge dual of \lstinline{self} with respect to the scalar product making $\lstinline{e}_0$, \dots, $\lstinline{e}_{\lstinline{n}}$ an orthonormal basis.

In order to endow the exterior algebra \lstinline{E} with a structure of Lie algebra, we construct a differential operator:
\begin{lstlisting}
sage: d = E.coboundary(str_eq)
\end{lstlisting}
where \lstinline{str_eq} is a dictionary whose keys are in $\{0,\dots,n\}\times\{0,\dots,n\}$ and whose values are (coerced into) $1$-forms in $\wedge^1 R\langle e_0, \dots, e_{n}\rangle$. More precisely, \lstinline!str_eq = { (J,K) : E.gens()[M] }! yields that $d (e_{M}) = e_{J} \wedge e_{K}$, the others being zero.

We can construct further structures, {\itshape e.g.} complex structures, scalar products, by lifting the corresponding structures to the exterior algebra with \lstinline{E.lift_morphism(mat)}, respectively \lstinline{E.lifted_bilinear_form(mat)} where \lstinline{mat} is the associated matrix with respect to the basis $(\lstinline{e}_0, \dots, \lstinline{e}_{\lstinline{n}})$.

\section{Computing cohomological invariants of nilmanifolds}
We start by computing the de Rham cohomology of a $6$-dimensional nilmanifold associated to the Lie algebra
\begin{eqnarray*}
\mathfrak{g}_{3.1}\oplus3\mathfrak{g}_{1} &:=& \left. \mathrm{span}\,\{e_0,\ldots,e_5\} \middle\slash \left\langle [e_2,e_3]=e_0 \right\rangle \right. .
\end{eqnarray*}
This is called $\mathfrak{g}_{3.1}\oplus3\mathfrak{g}_{1}$ in the notation by \cite{bochner}, and corresponds to $\mathfrak{h}_{8}:=(0, 0, 0, 0, 0, 12)$ in \cite{salamon}.
Since the structure equations have rational coefficients, by the Mal'tsev theorem \cite{malcev} the associated connected simply-connected Lie group admits a lattice, and then $\mathfrak{g}_{3.1}\oplus3\mathfrak{g}_{1}$ is the Lie algebra of a nilmanifold $X$. By \cite[Theorem 1]{nomizu}, the de Rham cohomology of $X$ is isomorphic to the cohomology of the complex $(\wedge^\bullet, d)$ where $d\alpha(\_,\_):=-\alpha([\_,\_])$ on $1$-forms and then extended by Leibniz rule.

We realize the Lie algebra:
\begin{lstlisting}
sage: E = ExteriorAlgebra(SR, 'e', 6)
sage: str_eq = {(2,3):-E.gens()[0], }
sage: d = E.coboundary(str_eq)
\end{lstlisting}
and check it:
\begin{lstlisting}
sage: [d(b) for b in E.gens()]
[-e2^e3, 0, 0, 0, 0, 0]
sage: all([d(d(b)) == 0 for b in E.gens()])
True
\end{lstlisting}
We construct the corresponding chain complex:
\begin{lstlisting}
sage: cplx_dR = d.chain_complex(); cplx_dR
Chain complex with at most 7 nonzero terms over Symbolic Ring
\end{lstlisting}
The method \lstinline{ascii_art(cplx_dR)} gives a pictorial representation of the complex.
We compute the cohomology of the Lie algebra by \lstinline{homology()}:
\begin{lstlisting}
sage: cplx_dR.homology()
{0: Vector space of dimension 1 over Symbolic Ring, \
    1: Vector space of dimension 5 over Symbolic Ring, \
    2: Vector space of dimension 11 over Symbolic Ring, \
    3: Vector space of dimension 14 over Symbolic Ring, \
    4: Vector space of dimension 11 over Symbolic Ring, \
    5: Vector space of dimension 5 over Symbolic Ring, \
    6: Vector space of dimension 1 over Symbolic Ring}
\end{lstlisting}
This says that
\begin{eqnarray*}
\dim_{\mathbb{R}} H^0_{dR}(X;\mathbb{R}) &=& 1 , \\[5pt]
\dim_{\mathbb{R}} H^1_{dR}(X;\mathbb{R}) &=& 5 , \\[5pt]
\dim_{\mathbb{R}} H^2_{dR}(X;\mathbb{R}) &=& 11 , \\[5pt]
\dim_{\mathbb{R}} H^3_{dR}(X;\mathbb{R}) &=& 14 , \\[5pt]
\dim_{\mathbb{R}} H^4_{dR}(X;\mathbb{R}) &=& 11 , \\[5pt]
\dim_{\mathbb{R}} H^5_{dR}(X;\mathbb{R}) &=& 5 , \\[5pt]
\dim_{\mathbb{R}} H^6_{dR}(X;\mathbb{R}) &=& 1 .
\end{eqnarray*}

We now compute the Poincaré polynomial of all the $6$-dimensional nilpotent Lie algebras, summarized in the dictionary \lstinline{alg_nilp_6} in Appendix \ref{app:algnilp6}. Recall that the Poincaré polynomial is the polynomial $\sum_k b_k(X) \cdot x^k$, where $b_k(X):=\dim_{\mathbb{R}}H^k(X;\mathbb{R})$.
\begin{lstlisting}
sage: E = ExteriorAlgebra(SR,'e',6)
sage: for algebra in alg_nilp_6.keys():
          d = E.coboundary(alg_nilp_6[algebra])
          H = d.chain_complex().homology()
          Poincare_poly = sum([H[j].dimension() * x^j \
              for j in range(len(H))])
          print algebra, "\n\t", Poincare_poly
\mathfrak{g}_{6.N18}^{-1} :
	x^6 + 2*x^5 + 4*x^4 + 6*x^3 + 4*x^2 + 2*x + 1
6\mathfrak{g}_{1} :
	x^6 + 6*x^5 + 15*x^4 + 20*x^3 + 15*x^2 + 6*x + 1
\mathfrak{g}_{3.1}\oplus 3g_{1} :
	x^6 + 5*x^5 + 11*x^4 + 14*x^3 + 11*x^2 + 5*x + 1
\mathfrak{g}_{6.N9} :
	x^6 + 3*x^5 + 5*x^4 + 6*x^3 + 5*x^2 + 3*x + 1
\mathfrak{g}_{6.N7} :
	x^6 + 3*x^5 + 6*x^4 + 8*x^3 + 6*x^2 + 3*x + 1
\mathfrak{g}_{6.N5} :
	x^6 + 4*x^5 + 8*x^4 + 10*x^3 + 8*x^2 + 4*x + 1
2\mathfrak{g}_{3.1} :
	x^6 + 4*x^5 + 8*x^4 + 10*x^3 + 8*x^2 + 4*x + 1
\mathfrak{g}_{6.N19} :
	x^6 + 2*x^5 + 3*x^4 + 4*x^3 + 3*x^2 + 2*x + 1
\mathfrak{g}_{5.2}\oplus g_{1} :
	x^6 + 3*x^5 + 5*x^4 + 6*x^3 + 5*x^2 + 3*x + 1
\mathfrak{g}_{6.N2} :
	x^6 + 2*x^5 + 3*x^4 + 4*x^3 + 3*x^2 + 2*x + 1
\mathfrak{g}_{5.5}\oplus g_{1} :
	x^6 + 4*x^5 + 7*x^4 + 8*x^3 + 7*x^2 + 4*x + 1
\mathfrak{g}_{6.N6} :
	x^6 + 3*x^5 + 6*x^4 + 8*x^3 + 6*x^2 + 3*x + 1
\mathfrak{g}_{6.N17} :
	x^6 + 3*x^5 + 5*x^4 + 6*x^3 + 5*x^2 + 3*x + 1
\mathfrak{g}_{5.1}\oplus g_{1} :
	x^6 + 4*x^5 + 9*x^4 + 12*x^3 + 9*x^2 + 4*x + 1
\mathfrak{g}_{6.N11} :
	x^6 + 2*x^5 + 4*x^4 + 6*x^3 + 4*x^2 + 2*x + 1
\mathfrak{g}_{6.N15} :
	x^6 + 3*x^5 + 5*x^4 + 6*x^3 + 5*x^2 + 3*x + 1
\mathfrak{g}_{6.N8} :
	x^6 + 3*x^5 + 5*x^4 + 6*x^3 + 5*x^2 + 3*x + 1
\mathfrak{g}_{4.1}\oplus 2g_{1} :
	x^6 + 4*x^5 + 7*x^4 + 8*x^3 + 7*x^2 + 4*x + 1
\mathfrak{g}_{6.N4} :
	x^6 + 4*x^5 + 8*x^4 + 10*x^3 + 8*x^2 + 4*x + 1
\mathfrak{g}_{5.6}\oplus g_{1} :
	x^6 + 3*x^5 + 5*x^4 + 6*x^3 + 5*x^2 + 3*x + 1
\mathfrak{g}_{6.N20} :
	x^6 + 2*x^5 + 3*x^4 + 4*x^3 + 3*x^2 + 2*x + 1
\mathfrak{g}_{6.N18}^{1} :
	x^6 + 2*x^5 + 4*x^4 + 6*x^3 + 4*x^2 + 2*x + 1
\mathfrak{g}_{6.N3} :
	x^6 + 3*x^5 + 8*x^4 + 12*x^3 + 8*x^2 + 3*x + 1
\mathfrak{g}_{6.N1} :
	x^6 + 3*x^5 + 6*x^4 + 8*x^3 + 6*x^2 + 3*x + 1
\mathfrak{g}_{6.N10} :
	x^6 + 3*x^5 + 5*x^4 + 6*x^3 + 5*x^2 + 3*x + 1
\mathfrak{g}_{6.N16} :
	x^6 + 3*x^5 + 4*x^4 + 4*x^3 + 4*x^2 + 3*x + 1
\end{lstlisting}

We consider now the following left-invariant complex structure on $X$, equivalently, linear complex structure on $\mathfrak{g}_{3.1}\oplus3\mathfrak{g}_{1}$:
$$ J e_0 := e_1, \quad J e_2 := e_3, \qquad J e_4 := e_5, $$
and consequently
$$ J e_1 := -e_0, \quad J e_3 := -e_2, \qquad J e_5 := -e_4. $$
Recall that, in terms of the dual vector space $(\mathfrak{g}_{3.1}\oplus3\mathfrak{g}_{1})^\vee$ with dual basis $(e^j)_j$, we set $J\alpha(\_):=\alpha(J^{-1}\_)$; then we have:
$$ J e^0 := e^1, \quad J e^2 := e^3, \qquad J e^4 := e^5. $$
We define such a structure and extend it to the exterior algebra $\wedge^\bullet(\mathfrak{g}_{3.1}\oplus3\mathfrak{g}_{1})^\vee$:
\begin{lstlisting}
sage: mat_J = matrix(6,6,[[ 0,-1, 0, 0, 0, 0], \
                          [ 1, 0, 0, 0, 0, 0], \
                          [ 0, 0, 0,-1, 0, 0], \
                          [ 0, 0, 1, 0, 0, 0], \
                          [ 0, 0, 0, 0, 0,-1], \
                          [ 0, 0, 0, 0, 1, 0]])
sage: J = E.lift_morphism(mat_J)
sage: for b in E.gens():
          print b, "|-->", J(b)
e0 |--> e1
e1 |--> -e0
e2 |--> e3
e3 |--> -e2
e4 |--> e5
e5 |--> -e4
\end{lstlisting}
Up to $\mathbb{C}$-linear extension to $\wedge^{\bullet}(\mathfrak{g}_{3.1}\oplus3\mathfrak{g}_{1})^\vee \otimes \mathbb{C}$, we notice that the forms $\varphi^0:=e^0-\sqrt{-1}J(e^0)$, $\varphi^1:=e^2-\sqrt{-1}J(e^2)$, $\varphi^2:=e^4-\sqrt{-1}J(e^4)$ yield a basis for the $\sqrt{-1}$-eigenspace of $J$:
\begin{lstlisting}
sage: varphi = [None, None, None]
sage: for j in range(3):
          varphi[j] = E.gens()[2*j] - I * J(E.gens()[2*j])
          print varphi[j], "|-->", J(varphi[j])
e0 - I*e1 |--> I*e0 + e1
e2 - I*e3 |--> I*e2 + e3
e4 - I*e5 |--> I*e4 + e5
\end{lstlisting}
We introduce the following function to compute the conjugate of a form:
\begin{lstlisting}
def bar(form):
    """
    CliffordAlgebraElement -> CliffordAlgebraElement

    It return the conjugate of form.
    
    sage: E = ExteriorAlgebra(SR,'e',4)
    sage: bar(I*E.gens()[1]*E.gens()[2]+E.gens()[0])
    -I*e1^e2 + e0
    """
    return sum([form.interior_product(b).constant_coefficient()\
        .conjugate() * b for b in E.basis()])
\end{lstlisting}
The forms $\bar\varphi^0=e^1-\sqrt{-1}J(e^1)$, $\bar\varphi^1=e^3-\sqrt{-1}J(e^3)$, $\bar\varphi^2=e^5-\sqrt{-1}J(e^5)$ yield a basis for the $(-\sqrt{-1})$-eigenspace of $J$:
\begin{lstlisting}
sage: for j in range(3):
          print bar(varphi[j]), "|-->", J(bar(varphi[j]))
e0 + I*e1 |--> -I*e0 + e1
e2 + I*e3 |--> -I*e2 + e3
e4 + I*e5 |--> -I*e4 + e5
\end{lstlisting}
The differentials of $\varphi^j$, $\bar\varphi^j$ are:
\begin{lstlisting}
sage: for j in range(3):
          print "d(", varphi[j], ") =", d(varphi[j])
          print "d(", bar(varphi[j]), ") =", d(bar(varphi[j]))
d( e0 - I*e1 ) = -e2^e3
d( e0 + I*e1 ) = -e2^e3
d( e2 - I*e3 ) = 0
d( e2 + I*e3 ) = 0
d( e4 - I*e5 ) = 0
d( e4 + I*e5 ) = 0
\end{lstlisting}
That, in the basis $(\varphi^j,\bar\varphi^j)_j$, we have
$$ d\varphi^0=\frac{\sqrt{-1}}{2}\varphi^1\wedge\bar\varphi^1, \quad d\varphi^1=0,\quad d\varphi^2=0, $$
$$ d\bar\varphi^0=\frac{\sqrt{-1}}{2}\varphi^1\wedge\bar\varphi^1, \quad d\bar\varphi^1=0,\quad d\bar\varphi^2=0. $$

We realize the Dolbeault complex $(\wedge^{\bullet}(\mathfrak{g}_{3.1}\oplus3\mathfrak{g}_{1})^\vee, \overline\partial)$.
\begin{lstlisting}
sage: E.<varphi0, varphi1, varphi2, \
          barvarphi0, barvarphi1, barvarphi2> = ExteriorAlgebra(SR)
sage: str_eq = {(1,4): I/2 * varphi0 + I/2 * barvarphi0}
sage: delbar = E.coboundary(str_eq)
sage: [delbar(b) for b in E.gens()]
sage: all([delbar(delbar(b)) == 0 for b in E.gens()])
[1/2*I*varphi1^barvarphi1, 0, 0, 1/2*I*varphi1^barvarphi1, 0, 0]
True
\end{lstlisting}
And we compute the Dolbeault cohomology:
\begin{lstlisting}
sage: HDol = delbar.chain_complex().homology() ; HDol
{0: Vector space of dimension 1 over Symbolic Ring, \
    1: Vector space of dimension 5 over Symbolic Ring, \
    2: Vector space of dimension 11 over Symbolic Ring, \
    3: Vector space of dimension 14 over Symbolic Ring, \
    4: Vector space of dimension 11 over Symbolic Ring, \
    5: Vector space of dimension 5 over Symbolic Ring, \
    6: Vector space of dimension 1 over Symbolic Ring}
\end{lstlisting}
Since the dimension of the Dolbeault cohomology and the dimension of the de Rham cohomology coincide, this means that the Fr\"olicher spectral sequence of the nilmanifold $X$ associated to $\mathfrak{g}_{3.1}\oplus3\mathfrak{g}_{1}$ degenerates at the first page, according with \cite[Theorem 4.1]{ceballos-otal-ugarte-villacampa}.

As in \cite{ceballos-otal-ugarte-villacampa}, the Lie algebra $\mathfrak{g}_{6.N6}=\mathfrak{h}_{11}=(0,0,0,12,13,14+23)$ admits a continuous family of complex structures, characterized by the basis $(\varphi^0,\varphi^1,\varphi^2)$ of $(1,0)$-forms having differentials
$$
d \varphi^{0} = 0, \quad d\varphi^1=\varphi^0\wedge\bar\varphi^0,
$$
$$
d\varphi^2 = \varphi^0\wedge\varphi^1+B\cdot\varphi^0\wedge\bar\varphi^1+|B-1|\cdot\varphi^1\wedge\bar\varphi^0 ,
$$
varying $B\in\mathbb{R}\setminus\{0,1\}$.
With the code as before,
\begin{lstlisting}
sage: E.<varphi0, varphi1, varphi2, \
          barvarphi0, barvarphi1, barvarphi2> = ExteriorAlgebra(SR)
sage: _ = var('B')
sage: str_eq = {(0,3): varphi1 - barvarphi1, \
          (0,1): varphi2, (0,4): B * varphi2, \
          (1,3): abs(B - 1) * varphi2 - B.conjugate() * barvarphi2, \
          (3,4): barvarphi2, (0,4): -abs(B - 1) * barvarphi2}
delbar = E.coboundary(str_eq)
sage: HDol = delbar.chain_complex().homology() ; HDol
sage: HDol = delbar.chain_complex().homology()
sage: [H[1].dimension() for H in HDol.items()]
[1, 3, 5, 6, 5, 3, 1]
\end{lstlisting}
we just get the Dolbeault cohomology for the generic value of the parameter. In order to solve this problem, we proceed as follows. We construct by hand the matrices associated to the differential $d\colon \wedge^j \mathfrak{g}_{6.N6}^\vee \to \wedge^{j+1} \mathfrak{g}_{6.N6}^\vee$ with respect to the basis $(\varphi^j, \bar\varphi^j)_j$:
\begin{lstlisting}
sage: d_mat = {j: matrix(len(E.basis(j+1)), len(E.basis(j)), \
          [delbar(list(E.basis(j))[n]).\
          interior_product(list(E.basis(j+1))[m]).constant_coefficient() \
          for m in range(len(E.basis(j+1))) \
          for n in range(len(E.basis(j)))]) \
          for j in range(len(E.gens())+1)}
sage: d_mat[-1] = matrix(1,0,[])
\end{lstlisting}
The Dolbeault cohomology,
$$ H^j_{\overline\partial}(\mathfrak{g}_{6.N6})=\frac{\ker (\overline\partial \colon \wedge^j \mathfrak{g}_{6.N6}^\vee)}{\mathrm{rk}(\overline\partial \wedge^{j-1}\mathfrak{g}_{6.N6}^\vee \to \wedge^j \mathfrak{g}_{6.N6}^\vee)} , $$
can then be computed as
\begin{lstlisting}
sage: hDol = {}
sage: for j in range(len(E.gens())+1):
          V = VectorSpace(SR,len(E.basis(j)))
          Z = V.subspace(d_mat[j].transpose().kernel())
          B = Z.subspace([V([d_mat[j-1][h,k] \
              for h in range(d_mat[j-1].nrows())]) \
              for k in range(d_mat[j-1].ncols())])
          H = Z.quotient(B)
          hDol[j] = dim(H)
sage: hDol
{0: 1, 1: 3, 2: 5, 3: 6, 4: 5, 5: 3, 6: 1}
\end{lstlisting}
The last step is to modify the function \lstinline{rank} in order to take care of the parameters. A non efficient way is the following:
\begin{lstlisting}
def rank(M, var=None):
    """
    (Matrix, Expression) -> dic

    It computes the rank of a matrix M depending on parameters var.
    It returs a dictionary whose keys are the possible values
    of the rank obtained at the corresponding items.
    
    sage: rank(matrix(3, 2, [var('t'),0,0,1,0,0]), t)
    {0: [], 1: [[t == 0]], 2: [[t != 0]]}
    sage: rank(matrix(3, 2, [0,0,0,1,0,0]))
    1
    sage: rank(matrix(3, 2, [1,0,0,1,0,0]))
    2
    """
    if var == None:
        return M.rank()
    rk = {}
    for j in range(min([M.nrows(), M.ncols()])+1):
        rk[j] = solve([sum([b*b.conjugate() for b in M.minors(j)])!=0] \
            + [b==0 for b in M.minors(j+1)], var)
    return rk
\end{lstlisting}

We use the same approach to compute the Bott-Chern and Aeppli cohomology of a nilmanifold $X$ associated to $\mathfrak{g}_{3.1}\oplus3\mathfrak{g}_{1}$, equivalently, of $\mathfrak{g}_{3.1}\oplus3\mathfrak{g}_{1}$ with the linear complex structure determined by the complex structure equations $d\varphi^0=\frac{\sqrt{-1}}{2}\varphi^1\wedge\bar\varphi^1$, $d\varphi^1=0$, $d\varphi^2=0$.
We start by realizing the exterior algebra with the differentials $\partial$ and $\overline\partial$:
\begin{lstlisting}
sage: E.<varphi0, varphi1, varphi2, \
          barvarphi0, barvarphi1, barvarphi2> = ExteriorAlgebra(SR)
sage: str_eq_delbar = {(1,4): I/2 * varphi0}
sage: delbar = E.coboundary(str_eq_delbar)
sage: [delbar(b) for b in E.gens()]
[1/2*I*varphi1^barvarphi1, 0, 0, 0, 0, 0]
sage: str_eq_del = {(1,4): I/2 * barvarphi0}
sage: ddel = E.coboundary(str_eq_del)
sage: [ddel(b) for b in E.gens()]
[0, 0, 0, 1/2*I*varphi1^barvarphi1, 0, 0]
\end{lstlisting}
We construct the matrices associated to $\partial$ and $\overline\partial$:
\begin{lstlisting}
sage: delbar_mat = {j: matrix(len(E.basis(j+1)), len(E.basis(j)), \
          [delbar(list(E.basis(j))[n]).\
          interior_product(list(E.basis(j+1))[m]).constant_coefficient() \
          for m in range(len(E.basis(j+1))) \
          for n in range(len(E.basis(j)))]) \
          for j in range(len(E.gens())+2)}
sage: delbar_mat[-1] = matrix(1,0,[])
sage: delbar_mat[-2] = matrix(0,0,[])

sage: del_mat = {j: matrix(len(E.basis(j+1)), len(E.basis(j)), \
          [ddel(list(E.basis(j))[n]).\
          interior_product(list(E.basis(j+1))[m]).constant_coefficient() \
          for m in range(len(E.basis(j+1))) \
          for n in range(len(E.basis(j)))]) \
          for j in range(len(E.gens())+2)}
sage: del_mat[-1] = matrix(1,0,[])
sage: del_mat[-2] = matrix(0,0,[])
\end{lstlisting}
We compute the Dolbeault, Bott-Chern, Aeppli cohomology respectively:
\begin{lstlisting}
sage: hDol = {}
sage: for j in range(len(E.gens())+1):
          V = VectorSpace(SR,len(E.basis(j)))
          Z = V.subspace(d_mat[j].transpose().kernel())
          B = Z.subspace([V([d_mat[j-1][h,k] \
              for h in range(d_mat[j-1].nrows())]) \
              for k in range(d_mat[j-1].ncols())])
          H = Z.quotient(B)
          hDol[j] = dim(H)
sage: hDol
{0: 1, 1: 5, 2: 11, 3: 14, 4: 11, 5: 5, 6: 1}

sage: hBC = {}
sage: for j in range(len(E.gens())+1):
          V = VectorSpace(SR,len(E.basis(j)))
          Z = V.subspace(delbar_mat[j].transpose().kernel()).\
              intersection(V.subspace(del_mat[j].transpose().kernel()))
          B = Z.subspace([V([(del_mat[j-1]*delbar_mat[j-2])[h,k] \
              for h in range((del_mat[j-1]*delbar_mat[j-2]).nrows())]) \
              for k in range((del_mat[j-1]*delbar_mat[j-2]).ncols())])
          H = Z.quotient(B)
          hBC[j] = dim(H)
sage: hBC
{0: 1, 1: 4, 2: 10, 3: 16, 4: 14, 5: 6, 6: 1}

sage: hA = {}
sage: for j in range(len(E.gens())+1):
          V = VectorSpace(SR,len(E.basis(j)))
          Z = V.subspace((del_mat[j+1]*delbar_mat[j]).transpose().kernel())
          B = Z.subspace([V([delbar_mat[j-1][h,k] \
              for h in range(delbar_mat[j-1].nrows())]) \
              for k in range(delbar_mat[j-1].ncols())] \
              + [V([del_mat[j-1][h,k] \
              for h in range(del_mat[j-1].nrows())]) \
              for k in range(del_mat[j-1].ncols())])
          H = Z.quotient(B)
          hA[j] = dim(H)
sage: hA
{0: 1, 1: 6, 2: 14, 3: 16, 4: 10, 5: 4, 6: 1}
\end{lstlisting}
The results can be found in \cite{angella-franzini-rossi, latorre-ugarte-villacampa}.

\section{Classifying lcs structures on Lie algebras}

In \cite{angella-bazzoni-parton}, we classify the locally conformally symplectic structures on $4$-dimensional Lie algebras up to linear equivalence. In this section, we perform the classification on one specific Lie algebra by using SageMath.

We consider the Lie algebra
$$ \mathfrak{r}_4 = (14+24,24+34,34,0) , $$
in Salamon's notation \cite{salamon}. That is, $\mathfrak{r}_4$ is the Lie algebra with basis $(e_0,e_1,e_2,e_3)$ and with non-zero brackets determined by $[e_0,e_3]=-e_1$, $[e_1,e_3]=-e_1-e_2$, $[e_2,e_3]=-e_2-e_3$. (Note that, to be coherent with the notation above, we use indices starting from $0$, while Salamon uses indices starting from $1$.) Equivalently, the dual $\mathfrak{r}_4^\vee$ has a basis $(e^0, e^1, e^2, e^3)$ such that
\begin{eqnarray*}
de^0=e^0\wedge e^3+e^1\wedge e^3, &\quad& de^1=e^1\wedge e^3+e^2\wedge e^3,\\[5pt]
de^2=e^2\wedge e^3, &\quad& de^3=0 .
\end{eqnarray*}

We create the differential complex $(\wedge^\bullet\mathfrak{r}_4^\vee, d)$ by
\begin{lstlisting}
sage: E = ExteriorAlgebra(SR, 'e', 4)
sage: str_eq = {
          (0,3) : E.gens()[0],
          (1,3) : E.gens()[0] + E.gens()[1],
          (2,3) : E.gens()[1] + E.gens()[2],
      }
sage: d = E.coboundary(str_eq)
\end{lstlisting}
We print the string of structure equations
\begin{lstlisting}
sage: print([d(b) for b in E.gens()])
[e0^e3 + e1^e3, e1^e3 + e2^e3, e2^e3, 0]
\end{lstlisting}
We check the Jacobi identity $d^2=0$:
\begin{lstlisting}
sage: all([(d*d)(b)==0 for b in E.gens()])
True
\end{lstlisting}

We want to decide which $1$-form is a possible Lee form, and which $2$-forms are lcs structures. So we construct a generic $1$-form, $\vartheta=\sum_{j=0}^3 \vartheta_j\cdot e^j$ where $\vartheta_j\in\mathbb{R}$:
\begin{lstlisting}
sage: thetacoeff = var(["theta%d" % i \
          for i in range(len(E.basis(1)))])
sage: theta = sum([thetacoeff[j] * E.gens()[j] \
          for j in range(len(E.basis(1)))])
\end{lstlisting}
We compute
\begin{lstlisting}
sage: d(theta)
theta0*e0^e3 + (theta0 + theta1)*e1^e3 + (theta1 + theta2)*e2^e3
\end{lstlisting}
We solve the condition $d(\vartheta)=0$ by
\begin{lstlisting}
sage: dtheta0 = solve([d(theta).interior_product(c).\
          constant_coefficient()==0 for c in E.basis(2)], \
          thetacoeff, solution_dict=True)
sage: dtheta0
[{theta0: 0, theta1: 0, theta2: 0, theta3: r1}]
\end{lstlisting}
that is, the generic $d$-closed $1$-form is $\vartheta = \vartheta_3 \cdot e^3$ for $\vartheta_3\in\mathbb{R}$. The case $\vartheta_3=0$ corresponds to $\vartheta=0$, that is, to the symplectic case; so, we are interested in the case $\vartheta_3\neq0$.
To make easier the simplification of forms according to the rules in a dictionary, or in a list coming from a \lstinline{solve} command, we provide the following:
\begin{lstlisting}
def simplify_form(phi, dic = {}):
    """
    (CliffordAlgebraElement, dic) -> CliffordAlgebraElement
    (CliffordAlgebraElement, list) -> CliffordAlgebraElement

    It simplify the form phi according to the substitutions in
    the dictionary dic.
    
    sage: E = ExteriorAlgebra(SR,'e',4)
    sage: var("a0 a01 a012")
    sage: phi = a0 * E.gens()[0] + a01 * E.gens()[0] \
        * E.gens()[1] + a012 * E.gens()[0] * E.gens()[1] \
        * E.gens()[2]
    sage: dic = {
        a0 : 1,
        a01 : 0,
        a012 : -1
    }
    sage: simplify_form(phi, dic)
    -e0^e1^e2 + e0
    """
    return sum([phi.interior_product(c).constant_coefficient().\
        subs(dic).simplify_full() * c for c in E.basis()])
\end{lstlisting}
When \lstinline{dic} is not provided, it simply applies \lstinline{simplify_full()} to the coefficients.
So we can simply type:
\begin{lstlisting}
sage: theta = simplify_form(theta, dtheta0[0]) ; theta
r1*e3
\end{lstlisting}
We save the new variable coefficients of $\vartheta$:
\begin{lstlisting}
sage: thetacoeff = [el for el in dtheta0[0].values() if el != 0]
\end{lstlisting}

We construct now the generic $2$-form $\Omega=\sum_{j<k}\omega_{jk} \cdot e^j \wedge e^k$:
\begin{lstlisting}
sage: Omegacoeff = var(["omega%d%d" % (i,j) \
          for i in range(len(E.basis(1))) \
          for j in range(i+1,len(E.basis(1)))])
sage: Omega = sum([Omegacoeff[j] * list(E.basis(2))[j] \
          for j in range(len(E.basis(2)))])
\end{lstlisting}
We compute the twisted differential $d_\vartheta(\Omega)=d(\Omega)-\vartheta\wedge\Omega$:
\begin{lstlisting}
sage: latex(d(Omega) - theta * Omega)
\end{lstlisting}
getting
\begin{multline*}
\left( -\omega_{01} r_{1} - 2 \, \omega_{01} \right)  e_{0} \wedge e_{1} \wedge e_{3} + \left( -\omega_{02} r_{1} - \omega_{01} - 2 \, \omega_{02} \right)  e_{0} \wedge e_{2} \wedge e_{3} \\[5pt]
+ \left( -\omega_{12} r_{1} - \omega_{02} - 2 \, \omega_{12} \right)  e_{1} \wedge e_{2} \wedge e_{3} .
\end{multline*}
We have now two possible cases.
\begin{itemize}
\item {\itshape Either $r_1\neq -2$}. In this case, we see that $\omega_{01}=\omega_{02}=\omega_{12}=0$. But if we compute
\begin{lstlisting}
sage: (simplify_form(Omega, \
          {omega01 : 0, omega02 : 0, omega12 : 0}))^2
\end{lstlisting}
we get \lstinline{0}, that is, $\Omega$ is degenerate in this case.
\item {\itshape Or $r_1=-2$}, that is, $\vartheta = -2 \cdot e^3$. In this case, we get the generic lcs structure:
\begin{lstlisting}
sage: theta = simplify_form(theta, {thetacoeff[0] : -2})
sage: dthOmega0 = solve([(d(Omega) - theta * Omega)\
          .interior_product(c).constant_coefficient() == 0 \
          for c in E.basis(3)], Omegacoeff, \
          solution_dict=True)
sage: Omega = simplify_form(Omega, dthOmega0[0])
sage: print "Omega =", Omega
sage: print "such that", (Omega * Omega)\
         .interior_product(E.basis(4)[0, 1, 2, 3]), " != 0"
Omega = r5*e0^e3 + r4*e1^e2 + r3*e1^e3 + r2*e2^e3
such that 2*r4*r5  != 0
\end{lstlisting}
We save the new variable coefficients of $\Omega$:
\begin{lstlisting}
sage: Omegacoeff = [el for el in dthOmega0[0].values() \
          if el != 0]
\end{lstlisting}
\end{itemize}

At the end, we are reduced to the generic lcs (non-symplectic) structure
$$ \vartheta = -2 \cdot e^3, \qquad \Omega = r_5 \cdot e^0 \wedge e^3 + r_4 \cdot e^1 \wedge e^2 + r_3 \cdot e^1 \wedge e^3 + r_2 \cdot e^2 \wedge e^3 ,$$
where $r_2, r_3, r_4, r_5\in\mathbb{R}$ satisfy $r_4\cdot r_5\neq0$.
We want now to decide which of the above forms are equivalent, namely, obtained by means of automorphisms of the Lie algebra. We start by realizing a generic morphism $\Psi$ of the exterior algebra:
\begin{lstlisting}
sage: mat = matrix(len(E.gens()), len(E.gens()), \
          [var("a%d%d" % (i,j)) for i in range(len(E.gens())) \
          for j in range(len(E.gens()))])
sage: morphism = E.lift_morphism(mat)
\end{lstlisting}
The conditions for \lstinline{morphism} to be an automorphism of the Lie algebra are $d(\Psi(e^j))=\Psi(d(e^j))$ for any $j\in\{0,1,2,3\}$ and $\det\Psi\neq0$:
\begin{lstlisting}
sage: [(d(morphism(b)) - morphism(d(b))).interior_product(c) \
          for b in E.gens() for c in E.basis(2)]
sage: print
sage: mat.det()
[a03*a10 + a03*a11 - a00*a13 - a01*a13,
 a03*a20 + a03*a21 - a00*a23 - a01*a23,
 a03*a30 + a03*a31 - a00*a33 - a01*a33 + a00,
 a13*a20 + a13*a21 - a10*a23 - a11*a23,
 a13*a30 + a13*a31 - a10*a33 - a11*a33 + a00 + a10,
 a23*a30 + a23*a31 - a20*a33 - a21*a33 + a10 + a20,
 a03*a11 + a03*a12 - a01*a13 - a02*a13,
 a03*a21 + a03*a22 - a01*a23 - a02*a23,
 a03*a31 + a03*a32 - a01*a33 - a02*a33 + a01,
 a13*a21 + a13*a22 - a11*a23 - a12*a23,
 a13*a31 + a13*a32 - a11*a33 - a12*a33 + a01 + a11,
 a23*a31 + a23*a32 - a21*a33 - a22*a33 + a11 + a21,
 a03*a12 - a02*a13,
 a03*a22 - a02*a23,
 a03*a32 - a02*a33 + a02,
 a13*a22 - a12*a23,
 a13*a32 - a12*a33 + a02 + a12,
 a23*a32 - a22*a33 + a12 + a22,
 0,
 0,
 a03,
 0,
 a03 + a13,
 a13 + a23]

a03*a12*a21*a30 - a02*a13*a21*a30 - a03*a11*a22*a30 \
    + a01*a13*a22*a30 + a02*a11*a23*a30 - a01*a12*a23*a30 \
    - a03*a12*a20*a31 + a02*a13*a20*a31 + a03*a10*a22*a31 \
    - a00*a13*a22*a31 - a02*a10*a23*a31 + a00*a12*a23*a31 \
    + a03*a11*a20*a32 - a01*a13*a20*a32 - a03*a10*a21*a32 \
    + a00*a13*a21*a32 + a01*a10*a23*a32 - a00*a11*a23*a32 \
     - a02*a11*a20*a33 + a01*a12*a20*a33 + a02*a10*a21*a33 \
     - a00*a12*a21*a33 - a01*a10*a22*a33 + a00*a11*a22*a33
\end{lstlisting}
We consider now the matrix
$$
\lstinline{mat}_1 =
\left(\begin{array}{cccc}
\frac{1}{r_{5}} & 0 & 0 & 0 \\
0 & \frac{1}{r_{5}} & 0 & 0 \\
0 & 0 & \frac{1}{r_{5}} & 0 \\
0 & \frac{r_{2}}{r_{4}} & -\frac{r_{3}}{r_{4}} & 1
\end{array}\right) ,
$$
(note that $r_5\neq0$ by non-degeneracy) and the automorphism $\Psi_1$ associated to $\lstinline{mat}_1$ in the basis $(e^1,e^2,e^3,e^4)$:
\begin{lstlisting}
sage: mat1 = matrix(4, 4, [ 1/Omegacoeff[2], 0, 0, 0, \
          0, 1 / Omegacoeff[2], 0, 0, \
          0, 0, 1 / Omegacoeff[2], 0, \
          0, Omegacoeff[0] / Omegacoeff[1], -Omegacoeff[3] / Omegacoeff[1], 1 ])
sage: morphism1 = E.lift_morphism(mat1)
\end{lstlisting}
We compute how it transform $\vartheta$ and $\Omega$:
\begin{lstlisting}
sage: theta = morphism1(theta)
sage: Omega = morphism1(Omega)
sage: theta ; Omega
-2*e3
e0^e3 + r4/r5^2*e1^e2
\end{lstlisting}
We have thus the lcs structures
$$ \left\{
\begin{array}{rcll}
\vartheta &=& -2\cdot e^3 &\\[5pt]
\Omega &=& e^0\wedge e^3+\sigma \cdot e^1\wedge e^3,\quad \text{for }\sigma\neq0 \; .
\end{array}
\right. $$

We conclude by showing that the above lcs forms are not equivalent, up to automorphisms of the Lie algebra.
More precisely, we take $\Omega_1:=e^0\wedge e^3+\sigma_1 \cdot e^1\wedge e^3$ (with $\sigma_1\neq0$) and $\Omega_2:=e^0\wedge e^3+\sigma_2 \cdot e^1\wedge e^3$ (with $\sigma_2\neq0$). We have to show that, if $\sigma_1\neq\sigma_2$, then there is no automorphism of the Lie algebra preserving $\theta:=-2\cdot e^3$ and transforming $\Omega_1$ into $\Omega_2$. We consider a general linear map $\psi$ associated to the matrix $A=(a_{jk})_{j,k}$ with respect to the basis $(e^j)_j$, and we extend it to the exterior algebra $\lstinline{E}$ by the Leibniz formula. In the ring $\mathbb{Q}[a_{00},a_{01},\dots,a_{33},\sigma_1,\sigma_2]$, we consider the ideal generated by the conditions that: $\psi$ is actually a morphism of the Lie algebra, that is, it commutes with the differential; $\psi$ preserves $\theta$, and transforms $\Omega_1$ into $\Omega_2$. An automorphism as above corresponds to a zero of this ideal with the condition of being invertible, that is, $\det A\neq0$. To make computations easier, we compute the Gr\"obner basis of the ideal. In SageMath, we get this as follows:
\begin{lstlisting}
sage: var_Omega = var("sigma1 sigma2")
sage: var_mat = var(["a%d%d" % (i,j) \
          for i in range(len(E.gens())) \
          for j in range(len(E.gens()))])

sage: Anello = QQ[var_mat + var_Omega]

sage: theta = -2 * E.gens()[3]
sage: Omega1 = E.gens()[0] * E.gens()[3] \
          + sigma1 * E.gens()[1] * E.gens()[2]
sage: Omega2 = E.gens()[0] * E.gens()[3] \
          + sigma2 * E.gens()[1] * E.gens()[2]

sage: mat = matrix(len(E.gens()), len(E.gens()), list(var_mat))
sage: morphism = E.lift_morphism(mat)

sage: ideal0 =[(d(morphism(b)) - morphism(d(b))).interior_product(c) \
          for c in E.basis(2) for b in E.basis(1)]
sage: ideal1 = [(morphism(theta) - theta).interior_product(c) \
          for c in E.basis(1)]
sage: ideal1 += [(morphism(Omega1) - Omega2).interior_product(c) \
          for c in E.basis(2)]

sage: B = Anello.ideal(ideal0 + ideal1).groebner_basis()
\end{lstlisting}
By \lstinline{latex(B)}, we get:
\begin{multline*}
\left[a_{21}^{2} + a_{31} \sigma_{2} -  a_{20}, a_{32} \sigma_{2} + a_{21}, a_{00} - 1, a_{01}, a_{02}, a_{03}, a_{10} -  a_{21},\right.\\[5pt]
\left. a_{11} - 1, a_{12}, a_{13}, a_{22} - 1, a_{23}, a_{33} - 1, \sigma_{1} -  \sigma_{2}\right]
\end{multline*}
In particular, since $\sigma_1-\sigma_2$ belongs to the ideal, it follows that $\Omega_1$ and $\Omega_2$ are never equivalent for $\sigma_1\neq\sigma_2$.

Summarizing, we have proven the following result.
\begin{proposition}[{see \cite{angella-bazzoni-parton}}]
Consider the $4$-dimensional Lie algebra $\mathfrak{r}_4=(14+24,24+34,34,0)$, and the dual basis $(e^0,e^1,e^2,e^3)$. There are infinite non-equivalent locally conformally symplectic structures on $\mathfrak{r}_4$, parametrized by $\sigma\in\mathbb{R}\setminus\{0\}$:
$$ \left\{
\begin{array}{rcl}
\vartheta &=& -2\cdot e^3 \\[5pt]
\Omega &=& e^0\wedge e^3+\sigma \cdot e^1\wedge e^2,
\end{array}
\right. $$
up to automorphisms of the Lie algebra.
\end{proposition}

Finally, we answer to the question whether such linear lcs structures correspond to invariant structures on smooth compact manifolds; that is, whether the connected simply-connected Lie group naturally associated to $\mathfrak{r}_4$ admits compact quotients. By \cite[Lemma 6.2]{milnor}, a necessary condition is that $\mathfrak{r}_4$ should be {\em unimodular}, that is, $\mathrm{tr}\,\mathrm{ad}\,X=0$ for all $X\in\mathfrak{r}_4$. At the level of the dual Lie algebra $\mathfrak{r}_4^\vee$, this can be checked by the vanishing of $d\wedge^3\mathfrak{r}_4^\vee$, otherwise there would exist an exact volume form.
By
\begin{lstlisting}
sage: [d(b) for b in E.basis(3)]
\end{lstlisting}
we get
\begin{lstlisting}
[3*e0^e1^e2^e3, 0, 0, 0]
\end{lstlisting}
showing that $\mathfrak{r}_4$ is not unimodular. Therefore, it is not associated to a compact manifold.

\appendix

\section{Six-dimensional nilpotent Lie algebras}\label{app:algnilp6}

Six-dimensional nilpotent Lie algebras are classified into $34$ different classes, up to isomorphism, by V.V. Morozov \cite{morozov}, see also \cite{magnin}, \cite[Table 15]{bochner}, see also \cite[Section 3]{gong}.
We report here their structure equations, in the notation of \cite{bochner}, in accord with the code in Section \ref{sec:preliminaries-sage}.

\begin{lstlisting}
alg_nilp_6 = {"\\mathfrak{g}_{6.N2}" : 
              {(0,4):E.gens()[1], (0,5):E.gens()[2], (0,2):E.gens()[3], \
              (0,3):E.gens()[4], },
              "\\mathfrak{g}_{6.N19}" : 
              {(5,2):E.gens()[1], (0,4):E.gens()[1], \
              (0,5):E.gens()[2]-E.gens()[3], (0,2):E.gens()[3], \
              (0,3):E.gens()[4], },
              "\\mathfrak{g}_{6.N11}" :
              {(5,2):E.gens()[1], (0,5):E.gens()[2]-E.gens()[3], \
              (0,2):E.gens()[3], (0,3):E.gens()[4], },
              "\\mathfrak{g}_{6.N18}^{1}" :
              {(2,4):E.gens()[0]+E.gens()[1], (0,2):E.gens()[1], \
              (0,4):E.gens()[3], (4,3):E.gens()[5], (1,2):E.gens()[5], },
              "\\mathfrak{g}_{6.N18}^{-1}" :
              {(4,2):E.gens()[0]+E.gens()[1], (2,0):E.gens()[1], \
              (4,0):E.gens()[3], (4,3):E.gens()[5], (2,1):E.gens()[5], },
              "\\mathfrak{g}_{6.N20}" :
              {(5,3):E.gens()[1], (0,4):E.gens()[1], \
              (0,5):2*E.gens()[2], \
              (0,2):1/2*E.gens()[3], (0,3):E.gens()[4], \
              (5,2):1/2*E.gens()[4],},
              "\\mathfrak{g}_{6.N6}" :
              {(0,2):E.gens()[1], (5,3):1/2*E.gens()[1], \
              (0,5):E.gens()[2], \
              (0,3):1/2*E.gens()[4], },
              "\\mathfrak{g}_{6.N7}" :
              {(0,3):E.gens()[1], (2,0):E.gens()[3], (2,4):E.gens()[5], },
              "\\mathfrak{g}_{6.N1}" :
              {(0,3):E.gens()[1], (0,2):E.gens()[3], (0,4):E.gens()[5], },
              "\\mathfrak{g}_{6.N3}" :
              {(0,4):E.gens()[1], (0,2):E.gens()[3], (2,4):E.gens()[5], },
              "\\mathfrak{g}_{6.N17}" :
              {(0,4):E.gens()[3], (2,1):E.gens()[3], (0,5):E.gens()[4], \
              (0,2): E.gens()[5], },
              "\\mathfrak{g}_{6.N15}" :
              {(0,4):E.gens()[3], (2,1):E.gens()[3], (2,5):E.gens()[3], \
              (0,5):E.gens()[4], (0,2):E.gens()[5], },
              "\\mathfrak{g}_{5.6}\\oplus g_{1}" :
              {(0,4):E.gens()[3], (2,5):E.gens()[3], (0,5): E.gens()[4], \
              (0,2):E.gens()[5], },
              "\\mathfrak{g}_{5.2}\\oplus g_{1}" :
              {(0,4):E.gens()[3], (0,5):E.gens()[4], (0,2):E.gens()[5], },
              "\\mathfrak{g}_{6.N9}" :
              {(0,5):E.gens()[1], (0,4):E.gens()[3], (5,2):E.gens()[3], \
              (0,2):E.gens()[5], },
              "\\mathfrak{g}_{6.N8}" :
              {(4,3):E.gens()[1], (4,2):E.gens()[1], (0,4):E.gens()[2], \
              (0,2):E.gens()[5], },
              "\\mathfrak{g}_{6.N16}" :
              {(0,3):E.gens()[1], (2,5):E.gens()[1], (0,5):E.gens()[3], \
              (2,4):E.gens()[3], (0,4):E.gens()[5], },
              "\\mathfrak{g}_{6.N10}" :
              {(0,5):E.gens()[1], (2,4):1/2*E.gens()[1], \
              (0,4):E.gens()[3], \
              (5,2):1/2*E.gens()[3], (0,2):1/2*E.gens()[5], },
              "\\mathfrak{g}_{4.1}\\oplus 2g_{1}" :
              {(0,3):E.gens()[1], (0,2):E.gens()[3], },
              "\\mathfrak{g}_{5.5}\\oplus g_{1}" :
              {(0,4):E.gens()[3], (2,5):E.gens()[3], (0,2):E.gens()[5], },
              "\\mathfrak{g}_{6.N4}" :
              {(0,4):E.gens()[3], (2,1):E.gens()[3], (0,2):E.gens()[5], },
              "2\\mathfrak{g}_{3.1}" :
              {(0,4):E.gens()[1], (2,5):E.gens()[3], },
              "\\mathfrak{g}_{5.1}\\oplus g_{1}" :
              {(0,4):E.gens()[1], (0,2):E.gens()[3], },
              "\\mathfrak{g}_{6.N5}" :
              {(0,5):E.gens()[1], (2,4):E.gens()[1], (0,4):E.gens()[3], \
              (5,2):E.gens()[3], },
              "\\mathfrak{g}_{3.1}\\oplus 3g_{1}" :
              {(0,2):E.gens()[1], },
              "6\\mathfrak{g}_{1}" :
              {},
              }
\end{lstlisting}

\end{document}